\documentclass[11pt,reqno]{article}

\usepackage{booktabs} 
\usepackage{array} 
\usepackage{paralist} 
\usepackage{verbatim} 
\usepackage{amsthm}
\usepackage[table,xcdraw]{xcolor}
\usepackage{amsmath}\usepackage{latexsym, amsmath, amssymb, a4, epsfig, color}
\usepackage{blindtext}
\usepackage{graphicx}
\usepackage{subcaption}
\usepackage[T1]{fontenc}
\usepackage[export]{adjustbox}
\usepackage{wrapfig}
\usepackage{chngcntr}
\usepackage{apptools}
\usepackage{amssymb}
\usepackage{float,subcaption}
\usepackage{stackengine}

\AtAppendix{\counterwithin{definition}{subsection}}

\graphicspath{}

\newtheorem{theorem}{Theorem}[section]
\newtheorem{cor}[theorem]{Corollary}
\newtheorem{lemma}[theorem]{Lemma}
\newtheorem{prop}[theorem]{Proposition}

\newtheorem{remark}{Remark}

\DeclareMathOperator{\supp}{supp}

 \def\f{\frac}
\newcommand{\nm}{\noalign{\smallskip}}
\newcommand{\ds}{\displaystyle}
\newcommand{\pf}{\noindent {\sl Proof}. \ }
\newcommand{\p}{\partial}

\def\bea{\begin{eqnarray}}  \def\eea{\end{eqnarray}}

\def\beas{\begin{eqnarray*}} \def\eeas{\end{eqnarray*}}

\newcommand{\Acal}{\mathcal{A}}

\def\Bx{{\bf x}}
\def\By{{\bf y}}

\newcommand{\eps}{\varepsilon}
\newcommand{\Ga}{\alpha}
\newcommand{\Gb}{\beta}

\newcommand{\Gl}{\lambda}

\newcommand{\RR}{\mathbb{R}}

\newcommand{\NN}{\mathbb{N}}

\newcommand{\Om}{\Omega}

\newcommand{\beq}{\begin{equation}}
\newcommand{\eeq}{\end{equation}}
\newcommand{\be}{\begin{equation*}}
\newcommand{\ee}{\end{equation*}}

\numberwithin{equation}{section}
\numberwithin{figure}{section}

\begin{document}

\title{A mathematical and numerical framework for near-field optics\thanks{\footnotesize
The work of D. S. Choi is supported by the Korean Ministry of Science, ICT and Future Planning through NRF grant No. 2016R1A2B4014530.}}
\author{
Habib Ammari\thanks{\footnotesize Department of Mathematics, ETH Z{\"u}rich, R{\"a}mistrasse 101, CH-8092  Z{\"u}rich, Switzerland (habib.ammari@math.ethz.ch, sanghyeon.yu@sam.math.ethz.ch).} \and 
Doo Sung Choi\thanks{\footnotesize Department of Mathematical Sciences, Korea Advanced Institute of Science and Technology, Daejeon 34141, South Korea (7john@kaist.ac.kr).} \and
Sanghyeon Yu\footnotemark[2]}

\date{}
\maketitle

\begin{abstract}
This paper is concerned with the inverse problem of reconstructing small and local perturbations of a planar surface
using the field interaction between a known plasmonic particle and the planar surface. 
The aim is to perform a super-resolved reconstruction of these perturbations from  shifts in the plasmonic frequencies
of the particle-surface system. In order to analyze the interaction between the plasmonic particle and the planar surface,   a well chosen conformal mapping, which transforms the particle-surface system into a coated structure, is used.  Then the even Fourier coefficients of the transformed domain are related to the shifts in the plasmonic resonances of the particle-surface system.  A direct reconstruction of the perturbations of the planar surface is proposed. Its viability and limitations are documented by numerical examples.   
\end{abstract}

\noindent {\footnotesize Mathematics Subject Classification
(MSC2000): 35R30, 35C20.}\\


\section{Introduction}

In conventional optical imaging and spectroscopy, a sample is typically illuminated by a light source and the scattered  light is recorded by a detector. The formed images are diffraction-limited. The diffraction limit essentially means that visible light cannot image  nanomaterials. 

In near-field optics, by placing a probe to exploit the unique properties of metal nanostructures at optical frequencies to localize incident illumination and enhance the light-matter interaction on the sample, one breaks the resolution limit \cite{bao, nf1,schotland,nf2,nf3,nf4,lukas1, lukas2}. Typically, the probe is a resonant plasmonic nanoparticle. At resonant frequencies, a light wave incident on the plasmonic probe  gives rise to a greatly amplified electric field just outside the probe, which then affects the nearby sample \cite{lukas1,lukas2}. When the plasmonic probe scans the sample surface one can form an image with a resolution much smaller than the diffraction limit by measuring the shifts of the plasmonic resonance induced by local defects for different particle's positions.  The physical mechanism is based on the excitation of plasmonic modes at the same particular frequencies  as the nanoprobe  and their coupling with the evanescent light.  The
plasmonic particle interacts with the surface and propagates its near field information into the far-field \cite{hai,hai2,hai3, lukas1,lukas2} through the introduced shifts in its plasmonic resonances.

The plasmon resonant frequency is one of the most important characteristics of a plasmonic particle. It depends not only on the electromagnetic properties of the particle and its size and shape \cite{matias, matias2, kelly, tocho}, but also on the electromagnetic properties of the environment \cite{matias, kelly, link}. It is the last property which enables  sensing applications of plasmonic particles \cite{link}.

In this paper, we  first provide a mathematical and computational framework to elucidate physical mechanisms for going beyond the diffraction limit in near-field optics based on the excitation of surface plasmons. We  mathematically and numerically analyze the intriguing behavior of light under the influence of plasmonics which allows nano-sensing of samples by using a localized surface plasmonic mode at the nanoprobe. We consider a plasmonic nanoparticle placed near a locally perturbed planar surface. We propose a mathematical and numerical framework to quantitatively image the sample from  shifts in the plasmonic resonances due to the coupling between the nanoparticle and the perturbed planar surface. The key idea is to use a well-chosen conformal mapping and to express the far-field and the shifts in the plasmonic resonances in terms of the Fourier coefficients of the transformed domain.  We compare the reconstructed images to those obtained from 
the contracted generalized polarization tensors of the transformed domain, which are the gold standard in wave imaging of small particles \cite{ammari2013mathematical,gpt3,gpt2, part1}. 

It is worth emphasizing that the physical system in which our method is realized is a plasmonic nanoparticle in the vicinity of a planar surface. The physical systems for imaging sub-wavelength surface defects described in \cite{nf3,nf4} can be mapped onto such a configuration.  Moreover, in near-field optics, besides the method based on excitation of surface plasmon described in this paper, there are have been a variety of attempts to attain super-resolution imaging. These include resonant energy transfer between fluorescence molecules and metallic or dielectric particles and interfaces and the generation or conversion of an evanescent wave at a grating with a period finer than the operating wavelength. We refer the interested reader to, for instance,  \cite{kawata} for an overview of these super-resolution methods in optics.

Our results in this work extend those in  \cite{part1,part2}, where a two particle system is considered. In \cite{part1,part2}, the system is composed of a known plasmonic particle and a small object.  By varying the relative position of the  particles, it is shown that fine details of the shape of the small object can be reconstructed from the induced  shifts of the plasmonic resonant frequencies of the plasmonic particle.

The paper is organized as follows. In Section \ref{sec-prelim}, we provide basic results on layer potentials and then explain the concepts of plasmonic resonances and contracted generalized polarization tensors.
In Section \ref{sec-forward}, we consider the forward scattering problem of the incident field interaction with a system composed of  a plasmonic particle and a perturbed planar surface. We apply a M\"{o}bius transformation in order to transform the plasmonic particle and the half-plane into concentric disks. Then we derive the asymptotic expansions of the scattered field and the plasmonic resonances of the particle-surface system. 
In Section \ref{sec-inverse}, we consider the inverse problem of reconstructing the perturbations of the planar surface. This is done  by relating the the shifts  in the resonances of the particle-planar surface system to the Fourier coefficients of the transformed domain. 
In Section \ref{sec-numerics}, we provide numerical examples to justify our theoretical results and to illustrate the performances of the proposed Fourier-based reconstruction scheme. In particular, we compare the reconstructed images to those obtained from 
the contracted generalized polarization tensors of the transformed domain.

\section{Preliminary results} \label{sec-prelim}

\subsection{Layer potentials}
We recall some basic results from layer potential theory that are needed for subsequent analysis. We refer to  \cite{book2} for more details.  
We denote by $\Gamma(\Bx,\By)$ the fundamental solution of the Laplacian in $\mathbb{R}^2$, i.e., 
$$\Gamma(\Bx,\By) = \frac{1}{2\pi}\ln |\Bx-\By|.$$
Let $D$ be a bounded domain in $\mathbb{R}^2$  with $\mathcal{C}^{1,\eta}$ boundary for some $\eta>0$, and let $\nu(\Bx)$ be the outward normal for $\Bx \in \partial D$. 

The single-layer potential $\mathcal{S}_{D}$ associated with $D$ is defined by 
$$
\mathcal{S}_{D} [\varphi](\Bx) =\int_{\p D }  \Gamma(\Bx,\By)\varphi(\By) d\sigma(\By) ,   \quad \Bx \in \mathbb{R}^2, 
$$
and the Neumann-Poincar\'{e} (NP) operator $\mathcal{K}_{D}^*$ by
$$
\mathcal{K}_{D}^* [\varphi](\Bx) = \int_{\p D }  \frac{ \p \Gamma }{\p\nu(\Bx)} (\Bx,\By) \varphi(\By) d\sigma(\By) ,   \quad \Bx \in \partial D.
$$
The following jump relations hold:
\begin{align}
{\mathcal{S}_D[\varphi]}\big|_+ &= {\mathcal{S}_D[\varphi]}\big|_-,\label{eqn_jump_single1}
\\
\frac{\p\mathcal{S}_D[\varphi]}{\p\nu}\Big|_{\pm} &= (\pm\frac{1}{2} I +\mathcal{K}_D^*)[\varphi]. \label{eqn_jump_single2}
\end{align}
Here, the subscripts $+$ and $-$ indicate the limits from outside and inside $D$, respectively.

Let $H^{1/2}(\p D)$ be the usual Sobolev space and let $H^{-1/2}(\p D)$ be its dual space with respect to the duality pairing $(\cdot,\cdot)_{-\frac{1}{2},\frac{1}{2}}$. We denote by $H^{-1/2}_0(\p D)$ the collection of all $\varphi\in H^{-1/2}(\p D)$ such that $(\varphi,1)_{-\frac{1}{2},\frac{1}{2}}=0$.

The NP operator is bounded from $H^{-1/2}(\p D)$ into $H^{-1/2}(\p D)$.
Moreover, the operator $\lambda I - \mathcal{K}_D^*: L^2(\partial D)
\rightarrow L^2(\partial D)$ is invertible for any $|\lambda| > 1/2$ \cite{book2,kellogg}.
Although the NP operator is not self-adjoint on $L^2(\p D)$, it can be symmetrized on $H_0^{-1/2}(\p D)$ with a proper inner product \cite{kang1,matias}. In fact, 
let $\mathcal{H}^*(\p D)$ be the space $H^{-1/2}_0(\p D)$ equipped with the inner product $(\cdot,\cdot)_{\mathcal{H}^*(\p D)}$ defined by
\begin{equation} \label{add1}
(\varphi,\psi)_{\mathcal{H}^*(\p D)} =  -(\varphi,\mathcal{S}_D[\psi])_{-\frac{1}{2},\frac{1}{2}},
\end{equation}
for $\varphi,\psi\in H^{-1/2}(\p D)$.
Then using the Plemelj's symmetrization principle,
$$
\mathcal{S}_D\mathcal{K}_D^*=\mathcal{K}_D\mathcal{S}_D,
$$
it can be shown that the NP operator $\mathcal{K}_D^*$ is self-adjoint in $\mathcal{H}^*$ with the inner product $(\cdot,\cdot)_{\mathcal{H}^*(\p D)}$ \cite{kang1,shapiro}. Since $\mathcal{K}_D^*$ is also compact, 
it admits the following spectral decomposition
in $\mathcal{H}^*$,
\begin{equation} \label{spectral_decomposition_Kstar}
\mathcal{K}_D^* = \sum_{j=1}^\infty \lambda_j (\cdot,\varphi_j)_{\mathcal{H}^*} \; \varphi_j,
\end{equation}
where $\lambda_j$ are the eigenvalues of $\mathcal{K}_D^*$ and $\varphi_j$ are their associated eigenfunctions. Note that $|\lambda_j|<1/2$ for all $j \geq 1$.

\subsection{Plasmonic resonances of metallic particles}\label{subsec-plasmonic}

In this subsection, we are interested in the frequency regime where plasmonic resonances occur in free space.  In such a regime, the wavelength of the incident field is much greater than the size of the plasmonic particle. To simplify the analysis and better illustrate the main idea, we use the quasi-static approximation (by assuming the incident wavelength to be infinite) to model the interaction. More precisely, let $\Omega$ represent a plasmonic particle with permittivity $\eps_c$ embedded in the homogeneous space  $\mathbb{R}^2$ with permittivity $\varepsilon_0$. We consider the following transmission problem with given incident field $u^i$ which is harmonic in $\mathbb{R}^2$:
\begin{equation} \left\{
    \begin{array}{l}
    \ds    \nabla \cdot (\eps \nabla u) = 0 \qquad \text{ in }\; \mathbb{R}^2, \\
        \nm
        \ds
         (u - u^i)(\Bx) = O(|\Bx|^{-1}) \qquad \text{ as }\; |\Bx| \rightarrow \infty,
    \end{array} \right.
    \label{add2}
\end{equation}
where $\eps  = \eps_c \chi(\Omega) +  \eps_0\chi(\mathbb{R}^2 \backslash \overline{\Omega})$, and $\chi(\Omega)$ and $\chi(\mathbb{R}^2 \backslash \overline{\Omega})$ are  the characteristic functions of $\Omega$ and $\mathbb{R}^2 \backslash \overline{\Omega}$, respectively.  
The problem (\ref{add2}) describes the response of the plasmonic particle $\Omega$ to the illumination $u^i$ in the quasi-static limit. 

The total field $u$ outside of $\Omega$ can be represented by (see \cite{kang1})
\begin{equation}
    u = u^i+ \mathcal{S}_{\Omega} [\varphi]  \, ,
    \label{scattered}
\end{equation}
where the density $\varphi$ satisfies the boundary integral equation
\begin{equation} \label{phi_int_equation}
 (\lambda I - \mathcal{K}_\Omega^*)[\varphi] = \frac{\partial u^i}{\partial \nu}\Big|_{\p \Omega}.
\end{equation}
Here, the permittivity contrast $\lambda$ is given by
 \begin{equation} \label{deflambda} 
\lambda= \frac{\eps_c +\eps_0}{2(\eps_c -\eps_0)}.
\end{equation}  

Contrary to ordinary dielectric particles, the permittivity $\eps_c$ of the plasmonic particle has negative real parts. In fact, $\eps_c$ depends on the operating frequency $\omega$ and can be modeled by the following Drude's model
\begin{equation}\label{drude}
\eps_c = \eps_c(\omega) = \eps_0\Big(1 -\frac{\omega_p^2}{\omega(\omega+i\gamma)}\Big),
\end{equation}
where $\omega_p>0$ is called the plasma frequency and $\gamma>0$ is the damping parameter. Since the parameter $\gamma$ is typically very small,
$\eps_c(\omega)$ has a small imaginary part.

Now we discuss the plasmonic resonances.
By applying the spectral decomposition \eqref{spectral_decomposition_Kstar} of $\mathcal{K}_\Omega^*$ to the integral equation \eqref{phi_int_equation}, we obtain \be\label{varphi_spectral_decomposition}
u = u^i + \sum_{j=1}^\infty \frac{(\frac{\p u^i}{\p\nu},\varphi_j)_{\mathcal{H}^*(\p \Omega)}}{\lambda - \lambda_j} \mathcal{S}_\Omega[\varphi_j].
\ee
Recall that $\lambda_j$ are eigenvalues $\mathcal{K}_\Omega^*$ and they satisfy the condition that $|\lambda_j|<1/2$. 
For $\omega<\omega_p$, $\Re [\eps_c(\omega)]$ can take  negative values. Then it holds that $|\Re [\lambda(\omega)]|<1/2$. If there exists a frequency, say $\omega_j$, such that $\lambda(\omega_j)$ is close to an eigenvalue $\lambda_j$ of the NP operator and their difference is locally minimized. Provided that $(\frac{\p u^i}{\p\nu},\varphi_j)_{\mathcal{H}^*(\p \Omega)}\neq 0$ , the eigenmode $\varphi_j$ in \eqref{varphi_spectral_decomposition} will be fully excited and it dominates over other modes. As a result, the scattered field $u-u^i$ will show a pronounced peak at the frequency $\omega_j$. This phenomenon is called the plasmonic resonance and $\omega_j$ is called the plasmonic resonant frequency. We refer the reader to \cite{kang1,hyeonbae, matias} for the details. 

We also mention that, for a precise reconstruction, the plasmonic particles with low-loss. In other words, the damping parameter should be small. It is known that silver and gold have relatively low loss. If loss is high, then the location of resonant peaks can be quite different from the one corresponding to the eigenvalue $\lambda_j$.


\subsection{Contracted generalized polarization tensors}\label{subsec-CGPT}

In this subsection, we explain the concept of the generalized polarization tensors (GPTs) associated with a smooth bounded domain $D$ having a permittivity $\epsilon_D$. Let $u$ be the solution of (\ref{add2})  and let $\lambda_D$ be defined by (\ref{deflambda})  with $\eps_c$ replaced by $\eps_D$. 
The scattered field $u-u^i$ has the following far-field behavior  \cite[p. 77] {book2}
\be
    (u - u^i)(\Bx) =  \sum_{|\alpha|, |\beta| \leq 1 } \f{(-1)^{|\beta|}}{\alpha ! \beta!} \partial^\alpha u^i(0) M_{\alpha \beta}(\lambda_D, D) \partial^\beta \Gamma(\Bx), 
    \quad |\Bx| \rightarrow + \infty,
    \label{scattered2}
\ee
where $M_{\alpha \beta}(\lambda_D,D)$ is given by
$$M_{\alpha \beta}(\lambda_D, D) : =  \int_{\partial D} y^\beta (\lambda_D I - \mathcal{K}_D^*)^{-1} [\frac{\partial \Bx^\alpha}{\partial \nu}](\By)\, d\sigma(\By), \qquad \alpha, \beta \in \mathbb{N}^2.$$
Here, the coefficient $M_{\alpha \beta}(\lambda_D, D)$ is called the {\it generalized polarization tensor} \cite{book2}.

For 
$m \in \NN$, let $P_m(\Bx)$ be the complex-valued polynomial
\begin{equation}
P_m(\Bx) = r^m \cos m\theta +ir^m \sin m\theta, \label{eq:Pdef}
\end{equation}
with $\Bx = re^{i\theta}$ in the polar coordinates.
For $n$ and 
$m$ in $\NN$, we define the {\it contracted generalized polarization tensors }(CGPTs) to be the
following linear combinations of generalized polarization tensors using the homogeneous harmonic polynomials introduced in
\eqref{eq:Pdef}:
\begin{equation} \label{defCGPT}
\begin{array}{l}
\ds M^{cc}_{mn}(\lambda_D,D) = \int_{\partial D} \Re [ P_n] (\lambda_D I - \mathcal{K}_D^*)^{-1} [\frac{\partial \,\Re  [ P_m]}{\partial \nu}]\, d\sigma,  \\
\nm \ds
M^{cs}_{mn}(\lambda_D,D) = \int_{\partial D} \Im  [ P_n]  (\lambda_D I - \mathcal{K}_D^*)^{-1} [\frac{\partial \,\Re  [ P_m]}{\partial \nu}]\, d\sigma,\\
\nm
\ds
M^{sc}_{mn}(\lambda_D,D) = \int_{\partial D} \Re[P_n] (\lambda_D I - \mathcal{K}_D^*)^{-1} [\frac{\partial \,\Im [ P_m]}{\partial \nu}]\, d\sigma,\\
\nm \ds
M^{ss}_{mn}(\lambda_D,D) = \int_{\partial D} \Im[ P_n]  (\lambda_D I - \mathcal{K}_D^*)^{-1} [\frac{\partial \,\Im [ P_m]}{\partial \nu}]\, d\sigma. \end{array}
\end{equation}
We remark that the CGPTs defined above encode useful information about the shape of the particle $D$ and can be used for its reconstruction. See \cite{gpt1,book2,  GPTs, gpt2} for more details.

It is sometimes more convenient to use the complex version of the contracted GPTs defined by 
\beq \label{comGPT1}
\mathbb{N}_{nm}^{(1)}:=M_{nm}^{cc}-M_{nm}^{ss}+i(M_{nm}^{cs}+M_{nm}^{sc}),
\eeq
\beq \label{comGPT2}
\mathbb{N}_{nm}^{(2)}:=M_{nm}^{cc}+M_{nm}^{ss}-i(M_{nm}^{cs}-M_{nm}^{sc}),
\eeq
for $n,m\neq 0$.  Here, $M_{nm}^{cc},  M_{nm}^{ss},M_{nm}^{cs}$, and $M_{nm}^{sc}$ are defined by (\ref{defCGPT}).  

\subsection{Shape derivatives of CGPTs}

Let $D_{\delta}$ be a small perturbation of the domain $D$ defined as follows: there is $h\in \mathcal{C}(\p D)$ such that $\p D_{\delta}$ is given by
\beq\label{h}
\p D_{\delta} = \{\Bx + \delta h(\Bx)\nu(\Bx) : \Bx \in \p D \}. 
\eeq
Then we have the following formula for the shape perturbation of the complex GPTs \cite{ammari2013mathematical,gpt1}:
\be
\begin{array}{l}
\ds
\mathbb{N}_{nm}^{(2)}(\lambda_D,D_{\delta})-\mathbb{N}_{nm}^{(2)}(\lambda_D,D)\nonumber\\
\nm \quad \ds = \delta \left(\frac{\varepsilon_D }{\varepsilon_0}-1\right) \int_{\p D_{0}}h(\Bx)\left[ {\frac{\partial u_n}{\p \nu}}\Big|_{-}{\frac{\partial {v_{-m}}}{\p \nu}}\Big|_{-} + \frac{\varepsilon_0}{\varepsilon_D}{\frac{\partial u_n}{\p T}}\Big|_{-} {\frac{\partial {v_{-m}}}{\p T}}\Big|_{-}\right](\Bx)d\sigma(\Bx) + O(\delta^2),\label{add4}
\end{array} \ee
where $u_n$ and $v_{m}$ are respectively the solutions to the following transmission problems:

\beq\label{ueqn}
\begin{cases}
\ds \Delta u =0 \quad&\mbox{in } \RR^2\setminus \partial D,\\
\ds u\vert_{+} = u\vert_{-} &\mbox{on } \partial D,\\
\ds \varepsilon_{0} {\frac{\partial u}{\p \nu}}\Big|_{+} = \varepsilon_{D} {\frac{\partial u}{\p \nu}}\Big|_{-} &\mbox{on } \partial D,\\
\ds (u- H)(\Bx) =O(|\Bx|^{-1}) &\mbox{as }  |\Bx|\to \infty,
\end{cases}
\eeq
and
\beq\label{veqn}
\begin{cases}
\ds \Delta v =0 \quad&\mbox{in } \RR^2\setminus \partial D,\\
\ds \varepsilon_{D} v\vert_{+} = \varepsilon_{0} v\vert_{-} &\mbox{on } \partial D,\\
\ds {\frac{\partial v}{\p \nu}}\Big|_{-} =  {\frac{\partial v}{\p \nu}}\Big|_{-} &\mbox{on } \partial D,\\
\ds (v - F)(\Bx) =O(|\Bx|^{-1}) &\mbox{as }  |\Bx|\to \infty,
\end{cases}
\eeq
with $H(\Bx)=r^{|n|}e^{in\theta}$ and $F(\Bx)=r^{|m|}e^{im\theta}$. 

When $D$ is the unit disk, the shape derivative can be computed explicitly. 
One can easily see that the solutions $u_n$ and $v_{m}$ of \eqref{ueqn} and \eqref{veqn} are given by
\begin{align}\label{un}
u_n(\Bx)
&=\ds\begin{cases}
\frac{2 \varepsilon_0}{\varepsilon_0 + \varepsilon_D} r^{|n|}e^{in\theta} \quad &\mbox{if } |\Bx|=r<1,\\[3mm]
\Bigl(r^{|n|} + \frac{\varepsilon_0 - \varepsilon_D}{\varepsilon_0 + \varepsilon_D} r^{-|n|}\Bigr)e^{in\theta} \quad &\mbox{if } |\Bx|=r>1,
\end{cases}
\end{align}
and
\begin{align}\label{vm}
v_{m}(\Bx)
&=\ds\begin{cases}
\frac{2 \varepsilon_D}{\varepsilon_0 + \varepsilon_D} r^{|m|}e^{i m\theta} \quad &\mbox{if } |\Bx|=r<1,\\[3mm]
\Bigl(r^{|m|} + \frac{\varepsilon_0 - \varepsilon_D}{\varepsilon_0 + \varepsilon_D} r^{-|m|}\Bigr)e^{i m\theta} \quad &\mbox{if } |\Bx|=r>1.
\end{cases}
\end{align}
Then, in view of (\ref{add4}),  we obtain the following result.
\begin{lemma}\label{Fouriermoment}
For $n,m\neq 0$,
\begin{equation} \label{add5}
\mathbb{N}_{nm}^{(2)}(\lambda_{D},D_{\delta}) - \mathbb{N}_{nm}^{(2)}(\lambda_{D},D)
 = \delta \frac{2\pi(\varepsilon_D|nm| + \varepsilon_0 nm)}{(\varepsilon_D - \varepsilon_0)\lambda_{-+}^2} \hat{h}(-n+m) + O(\delta^2),
\end{equation}
where $\lambda_{D}$ is given as $\lambda_D = (\varepsilon_D+ \varepsilon_0 )/(2(\varepsilon_D-\varepsilon_0))$ and $\hat{h}(k)$ is the Fourier coefficient of $h(\theta):=h(\cos\theta,\sin\theta)$ defined by $$\hat{h}(k) = \frac{1}{2\pi}\int_0^{2\pi} h(\theta) e^{-ik\theta}d\theta.$$
\end{lemma}
For later use, we note that 
\beq\label{hconjugate}
\hat{h}(-k) = \overline{\hat{h}(k)}
\eeq
for all $k\in \NN$, which follows from the fact that $h$ is a real-valued function.

\section{The forward problem} \label{sec-forward}  
In this section, we consider a plasmonic nanoparticle placed close to a locally perturbed planar surface.  We  let $\mathbb{H}_{0}$ to be the (unperturbed) lower half-plane 
$$\mathbb{H}_{0} = \{ \Bx \in \mathbb{R}^2 : \Bx=(x_1,x_2),\quad x_2<0 \}$$
with $\partial \mathbb{H}_{0} = \{ \Bx \in \mathbb{R}^2 : \Bx=(x_1,0)  \}$ being its boundary. 

We define $\mathbb{H}_{\delta}$ to be a $\delta$-perturbation of $\mathbb{H}_{0}$, {\it i.e.}, we let $h_0 \in \mathcal{C}(\p \mathbb{H}_{0})$ and $\p \mathbb{H}_{\delta}$ be given by
$$\p \mathbb{H}_{\delta} = \{\Bx + \delta h_0(\Bx)\nu(\Bx) : \Bx \in \p \mathbb{H}_{0} \},$$
where
\beq\label{supph_0}
\supp(h_0)\subset [-R,R] \quad \mbox{for some } R>0.
\eeq
We assume that $\delta$ is a small positive parameter and $R$ is of order one.

We also define a plasmonic particle $\Om$:
$$\Om = \{ \Bx \in \mathbb{R}^2 : \lvert \Bx - (0,d ) \rvert < r_\Omega  \},$$
where $r_\Omega$ is the radius of the plasmonic particle and $d$ is the distance between the center of $\Om$ and $\mathbb{H}_{0}$. Assume that they are of order one, {\it i.e.}, $r_\Omega, d=O(1)$.

\subsection{Transmission problem in the perturbed half-space}

We consider the following transmission problem:
\begin{equation}
\label{add0} \left\{ \begin{array}{l}
\ds \nabla \cdot (\varepsilon \nabla u) =0 \quad \mbox{in } \RR^2\setminus (\partial \mathbb{H}_{\delta} \cup \partial \Om),\\
\nm
\ds u\vert_{+} = u\vert_{-} \quad \mbox{on } \partial \mathbb{H}_{\delta} \cup \partial \Om,\\
\nm
\ds \varepsilon_{+} {\frac{\partial u}{\p \nu}}\Big|_{+} = \varepsilon_{-} {\frac{\partial u}{\p \nu}}\Big|_{-} \quad \mbox{on } \partial \mathbb{H}_{\delta},\\
\nm
\ds \varepsilon_{+} {\frac{\partial u}{\p \nu}}\Big|_{+} = \varepsilon_{c} {\frac{\partial u}{\p \nu}}\Big|_{-} \quad \mbox{on } \partial \Om,\\
\nm
\ds (u- u^i)(\Bx) =O(|\Bx|^{-1}) \quad \mbox{as }  |\Bx|\to \infty,
\end{array} \right.
\end{equation}
where $\varepsilon = \varepsilon_- \chi(\mathbb{H}_{\delta}) + \varepsilon_+ \chi(\mathbb{R}^2 \setminus({\mathbb{H}_{\delta} \cup \Om})) + \varepsilon_c \chi(\Om)$ with $\varepsilon_-, \varepsilon_+ > 0$. Here, $u^i=-\mathbf{a}\cdot \mathbf{x}$ represents the potential of a uniform incident field in the direction $\mathbf{a}\in\mathbb{R}^2$.

\begin{figure*}
\centering
\begin{subfigure}{.5\textwidth}
  \centering
  \includegraphics[trim = 0 100 0 100, clip=true, width=\linewidth]{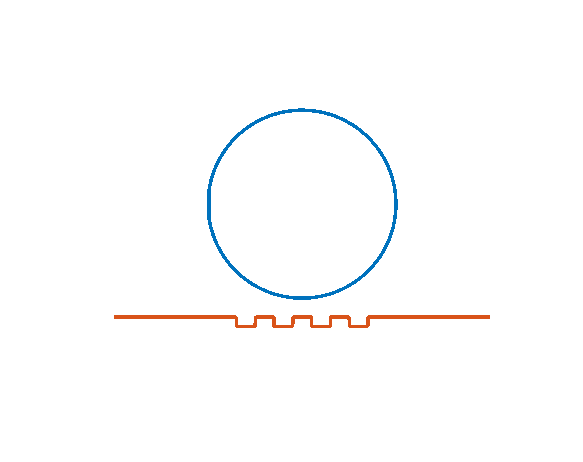}
  \caption{$h_0$ perturbation on the half-plane: $\p \mathbb{H}_{\delta}$}
  \label{a}
\end{subfigure}%
{\LARGE$\xrightarrow{\Phi}$}%
\begin{subfigure}{.5\textwidth}
  \centering
  \includegraphics[trim = 0 100 0 100, clip=true, width=\linewidth]{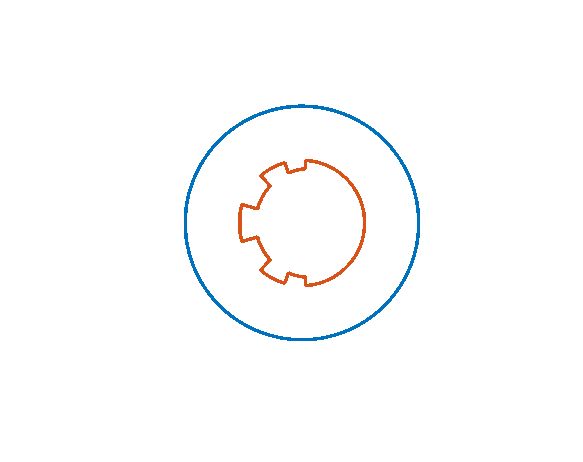}
  \caption{$h$ perturbation on the unit disk: $\p D_{1,\delta}$}
  \label{b}
\end{subfigure}
\caption{The plasmonic particle (left, blue) and the perturbed half-plane (left, red) are transformed into concentric perturbed disks (right).}
\label{transform}
\end{figure*}

To capture the interaction between the half-space and the plasmonic particle, we use a conformal mapping technique.
Define the conformal transformation $\Phi$, for $(x_1,x_2) \in \mathbb{R}^2$,  by
$$\Phi(z) = \frac{z+i\sqrt{d^2-r_\Omega^2}}{z-i\sqrt{d^2-r_\Omega^2}}, \quad z = x_1+ix_2.$$
Then $\mathbb{H}_{\delta}$ and $\mathbb{R}^2 \setminus\overline{\Om}$ are transformed into the domains \be
D_{1,\delta} := \Phi(\mathbb{H}_{\delta}) \mbox{ and } D_2 := \Phi(\mathbb{R}^2 \setminus\overline{\Om}) = \{ \lvert \zeta \rvert < e^s \},
\ee
where the parameter $s>0$ is given by
$$
 \sinh s = (\sqrt{(d/r_\Omega)^2-1}).
$$
Note that the domain $D_{1,\delta}$ becomes the unit disk when $\delta=0$, {\it i.e.},
$$D_{1,0} = \Phi(\mathbb{H}_{0}) = \{ \lvert \zeta \rvert < 1 \}.$$
Since $D_{1,\delta}$ is a $\delta$-perturbation of $D_{1,0}$, there is $h\in \mathcal{C}(\p D_{1,0})$ such that $\p D_{1,\delta}$ is given by
\beq\label{h2}
\p D_{1,\delta} = \{\Bx + \delta h(\Bx)\nu(\Bx) : \Bx \in \p D_{1,0} \} = \{e^{i\theta} + \delta h(\theta)e^{i\theta} : \theta \in [0,2\pi) \}.
\eeq
Figure \ref{transform} describes the transforming system of the whole domain.

For convenience, we denote $D_{1,\delta}$ by $D_1$. Since the mapping $\Phi$ is conformal, it can be shown that the transformed potential $\widetilde{u}(\zeta):= u (\Phi^{-1}(\zeta))$  satisfies the following transmission problem:
\be
\begin{cases}
\ds \nabla \cdot (\tilde{\varepsilon} \nabla \tilde{u}) =0 \quad&\mbox{in } \RR^2 \setminus (\partial{D_1} \cup \partial D_2),\\
\ds \tilde{u}\vert_{+} = \tilde{u}\vert_{-} &\mbox{on } \partial{D_1} \cup \partial D_2,\\
\ds {\varepsilon}_+ {\frac{\partial \tilde{u}}{\p \nu_1}}\Big|_{+} = {\varepsilon}_- {\frac{\partial \tilde{u}}{\p \nu_1}}\Big|_{-} &\mbox{on } \partial D_1,\\
\ds {\varepsilon}_{c} {\frac{\partial \tilde{u}}{\p \nu_2}}\Big|_{+} = {\varepsilon}_+ {\frac{\partial \tilde{u}}{\p \nu_2}}\Big|_{-} &\mbox{on } \partial{D_2},\\
\ds (\tilde{u}-\tilde{u}^i)(\zeta) =O(|\zeta-(1,0)|) &\mbox{as }  \zeta \to (1,0),
\end{cases}
\ee
where $\tilde{\varepsilon} = \varepsilon_- \chi(D_1) + \varepsilon_+ \chi(D_2 \setminus \overline{D}_1) + \varepsilon_c \chi(\mathbb{R}^2 \setminus D_2)$ and ${\frac{\partial}{\p \nu_i}}$ denotes the outward normal derivative with respect to $\p D_i$ for $i=1,2$.

The solution $\tilde{u}$ can be represented by
\begin{equation} \label{add3} \tilde{u} = \tilde{u}^i + \mathcal{S}_{D_1}[\phi_1] + \mathcal{S}_{D_2}[\phi_2],\end{equation}
where the densities $\phi_1$ and $\phi_2$ satisfy the system of integral equations
\begin{equation} \label{neumann1}
\left\{ \begin{array}{l}
\ds {\varepsilon}_+ {\frac{\partial}{\p \nu_1}}\Big( \tilde{u}^i + \mathcal{S}_{D_1}[\phi_1] + \mathcal{S}_{D_2}[\phi_2] \Big)\Big|_{+} = {\varepsilon}_- {\frac{\partial}{\p \nu_1}}\Big( \tilde{u}^i + \mathcal{S}_{D_1}[\phi_1] + \mathcal{S}_{D_2}[\phi_2] \Big)\Big|_{-} \quad \mbox{on } \partial D_1,\\
\nm
\ds
{\varepsilon}_{c} {\frac{\partial}{\p \nu_2}}\Big( \tilde{u}^i + \mathcal{S}_{D_1}[\phi_1] + \mathcal{S}_{D_2}[\phi_2] \Big)\Big|_{+} = {\varepsilon}_+ {\frac{\partial}{\p \nu_2}}\Big( \tilde{u}^i + \mathcal{S}_{D_1}[\phi_1] + \mathcal{S}_{D_2}[\phi_2] \Big)\Big|_{-} \quad \mbox{on } \partial D_2,
\end{array} \right.
\end{equation}
or equivalently, 
\begin{equation} \label{phicond1}
\left\{ \begin{array}{l}
\ds (\lambda_{-+}I - \mathcal{K}_{D_1}^*)[\phi_1] = {\frac{\partial}{\p \nu_1}}\Big( \tilde{u}^i + \mathcal{S}_{D_2}[\phi_2] \Big) \quad \mbox{on } \partial D_1,\\
\nm
\ds
(\lambda_{+c}I - \mathcal{K}_{D_2}^*)[\phi_2] = {\frac{\partial}{\p \nu_2}}\Big( \tilde{u}^i + \mathcal{S}_{D_1}[\phi_1] \Big) \quad \mbox{on } \partial D_2,
\end{array} \right.
\end{equation}
with
\beq\label{lambda}
\lambda_{-+} = \frac{\varepsilon_- + \varepsilon_+}{2(\varepsilon_- - \varepsilon_+)} \quad \mbox{and} \quad \lambda_{+c} = \frac{\varepsilon_+ + \varepsilon_{c}}{2(\varepsilon_+ - \varepsilon_{c})}.
\eeq

Since $|\lambda_{-+}| > \frac{1}{2}$ and $|\lambda_{+c}| < \frac{1}{2}$, the operator $(\lambda_{-+}I - \mathcal{K}_{D_1}^*)$ is invertible. Then we have from \eqref{phicond1} that
\beq\label{phi1}
\phi_1 = (\lambda_{-+}I - \mathcal{K}_{D_1}^*)^{-1}{\frac{\partial}{\p \nu_1}}\Big( \tilde{u}^i + \mathcal{S}_{D_2}[\phi_2] \Big).
\eeq
By substituting \eqref{phi1} into \eqref{phicond1}, we obtain the following result.
\begin{prop}
The density  $\phi_2$ on $\partial D_2$ satisfies the equation
\be
(\lambda_{+c}I - {\mathcal{A}} ) [\phi_2] = {\frac{\partial \tilde{u}_{D_1}}{\p \nu_2}} \quad \mbox{on } \partial D_2,
\ee
where the operator ${\mathcal{A}}$ is given by
\be
{\mathcal{A}} := \mathcal{K}_{D_2}^* + {\frac{\partial}{\p \nu_2}} \mathcal{S}_{D_1}(\lambda_{-+}I - \mathcal{K}_{D_1}^*)^{-1}{\frac{\partial \mathcal{S}_{D_2}[\cdot]}{\p \nu_1}},
\ee
and 
\be
\tilde{u}_{D_1} := \tilde{u}^i + \mathcal{S}_{D_1}(\lambda_{-+}I - \mathcal{K}_{D_1}^*)^{-1}\left[{\frac{\partial \tilde{u}^i}{\p \nu_1}} \right].
\ee
\end{prop}

It can be shown that the operator $\mathcal{A}$ is self-adjoint and compact \cite{part2}. So it has real eigenvalues and admits the following spectral decomposition:
\begin{equation}\label{spec_dec_A}
\mathcal{A} = \sum_{n\in \mathbb{Z}} \lambda_n \varphi_n\otimes \varphi_n,
\end{equation}
where $(\lambda_n,\varphi_n)$ are the eigenvalue-eigenfunction pairs of the operator $\mathcal{A}$. It can also be shown that the eigenvalues satisfy $|\lambda_n|\leq 1/2$.

Now we discuss how to measure the plasmonic resonance frequency or the eigenvalues $\lambda_n$. 
One can see that, using \eqref{spec_dec_A}, the polarizability $p_\Omega$ of the plasmonic particle $\Omega$ is given by
\begin{equation}\label{pol}
p_\Omega := \int_{\p \Omega} u^i \psi_{2} d\sigma
=  \sum_{n\in \mathbb{Z}} \frac{
\big( \tilde{u}^i,\varphi_{n}\big)_{L^2(\p D_2)}
\big( \varphi_n,{\frac{\partial \tilde{u}_{D_1}}{\p \nu_2}}\big)_{\mathcal{H}^*(\p D_2)}
}{\lambda_{+c}- \lambda_n},
\end{equation}
where $\psi_{2} = |\Phi'|^{-1} (\phi_2 \circ \Phi)$.
The absorption cross section, which can be measured from the far-field region, by the small particle $\Omega$ is proportional to the imaginary part of the polarizability $p_\Omega$ in the quasi-static limit. Recalling the Drude model \eqref{drude}, 
we see that
$\lambda_{+c}$ depends on the frequency  and satisfies $|\mbox{Re}\lambda_{+c}|\leq 1/2$ for $\omega<\omega_p$. So, the particle's resonance frequencies and the corresponding eigenvalues $\lambda_n$ can be measured as discussed in Subsection \ref{subsec-plasmonic}. 
From \eqref{pol}, it is clear that 
when we vary the frequency of the incident field, at the frequency $\omega_n$ such that $\lambda_{D_2}(\omega_n)= \lambda_n$ for some $n$ which satisfies the condition 
$$ 
\big( \tilde{u}^i,\varphi_{n}\big)_{L^2(\p D_2)}
\big( \varphi_n,{\frac{\partial \tilde{u}_{D_1}}{\p \nu_2}}\big)_{\mathcal{H}^*(\p D_2)}
\neq 0,
$$
the absorption cross section shows a sharp peak, which means that a plasmonic resonance is excited. In this way, we can measure the resonance frequency $\omega_n$ or the eigenvalue $\lambda_n$. In the next subsection, we analyze the asymptotic structure of the eigenvalue $\lambda_n$ for small $\delta$.

\subsection{Asymptotic expansion of the operator ${\mathcal{A}}$}
\label{subsec:asymp}
Here we compute an asymptotic expansion of the operator ${\mathcal{A}}:\mathcal{H}^*(\p D_2) \to \mathcal{H}^*(\p D_2)$ for small $\delta$, where $\mathcal{H}^*(\p D_2)$ is defined by (\ref{add1}) with $D$ replaced by $D_2$. Since $\p D_{1,0} $ is the unit circle, we use  the Fourier basis $\{ e^{in\theta} \}_{n\neq 0}$ as a basis of $\mathcal{H}^*(\p D_2)$. Let $(r,\theta)$ be the polar coordinates in the $\zeta$-plane, \textit{i}.\textit{e}., $\zeta = re^{i\theta}$. Then the following proposition holds.


\begin{prop}\label{Aexpansion}
We have the following asymptotic expansion of the operator $\mathcal{A}$ for small $\delta$ as follows: for $\zeta=e^{s+i\theta} \in \p D_2$,
\begin{align*}
\mathcal{A}[e^{in\theta}](\zeta)
=& -\frac{1}{4\lambda_{-+}}e^{-2|n|s}e^{in\theta} \\
&- \delta \sum_{m \in \mathbb{Z} \setminus \{0\}} \frac{(\varepsilon_-|nm| + \varepsilon_+nm)}{4 |n|(\varepsilon_- - \varepsilon_+)\lambda_{-+}^2} e^{-(|n|+|m|)s} \hat{h}(-n+m) e^{im\theta} + O(\delta^2)
\end{align*}
for $n\neq 0$.
\end{prop}
\pf
Following the same arguments as those in \cite{part2}, we obtain 
\begin{align*}
&\mathcal{A}[\varphi_n^c](\zeta) = \sum_{m =1}^\infty - \frac{1}{4\pi |n|}e^{-(|n|+m)s} [M_{nm}^{cc}(\lambda_{-+},D_{1,\delta})\varphi_m^c(\theta) + M_{nm}^{cs}(\lambda_{-+},D_{1,\delta})\varphi_m^s(\theta)],\\
&\mathcal{A}[\varphi_n^s](\zeta) = \sum_{m =1}^\infty - \frac{1}{4\pi |n|}e^{-(|n|+m)s} [M_{nm}^{sc}(\lambda_{-+},D_{1,\delta})\varphi_m^c(\theta) + M_{nm}^{ss}(\lambda_{-+},D_{1,\delta})\varphi_m^s(\theta)],
\end{align*}
for $\varphi_n^c(\theta) = \cos n\theta$, $\varphi_n^s(\theta) = \sin n\theta$, and $n\neq 0$.
Then, by using the linearity  of $\mathcal{A}$, it follows that
\begin{align}
\mathcal{A}[e^{in\theta}](\zeta)
=& \mathcal{A}[\varphi_n^c](\zeta) + i\mathcal{A}[\varphi_n^s](\zeta)\nonumber\\
=& \sum_{m \in \mathbb{Z} \setminus \{0\}}- \frac{1}{8\pi |n|} e^{-(|n|+|m|)s} \mathbb{N}_{nm}^{(2)}(\lambda_{-+},D_{1,\delta})e^{im\theta}.
\label{exact_A_mat}
\end{align}
We then decompose $\mathcal{A}[e^{in\theta}]$ as follows:
\begin{align*}
\mathcal{A}[e^{in\theta}](\zeta) &=
\sum_{m \in \mathbb{Z} \setminus \{0\}}- \frac{1}{8\pi |n|}e^{-(|n|+|m|)s} \mathbb{N}_{nm}^{(2)}(\lambda_{-+},D_{1,0})e^{im\theta} 
\\
& +\sum_{m \in \mathbb{Z} \setminus \{0\}}- \frac{1}{8\pi |n|}e^{-(|n|+|m|)s} [\mathbb{N}_{nm}^{(2)}(\lambda_{-+},D_{1,\delta})
-\mathbb{N}_{nm}^{(2)}(\lambda_{-+},D_{1,0})
]e^{im\theta} 
\\
&:=I + II.
\end{align*}
Note that $I$ is of order one but $II$ is of order $\delta$.
First, we compute $I$.
 Since $\p D_1$ and $\p D_2$ are circles,
\beq\label{vanishK}
\mathcal{K}_{D_1}^*[e^{in\theta}]= \mathcal{K}_{D_2}^*[e^{in\theta}]=0;
\eeq
see, for instance, \cite{ammari2013mathematical}. Moreover, from \cite{ammari2013spectral} it follows that
\begin{align}\label{normderivS}
\frac{\p}{\p r} \mathcal{S}_B[e^{in\theta}](w)
&=\ds\begin{cases}
-\frac{1}{2}\left(\frac{r}{r_0}\right)^{|n|}e^{in\theta} \quad &\mbox{if } |w|=r<r_0,\\
\frac{1}{2}\left(\frac{r_0}{r}\right)^{|n|}e^{in\theta} \quad&\mbox{if } |w|=r>r_0,
\end{cases}
\end{align}
where $B$ is the disk centered at the origin and with radius $r_0$.
By using \eqref{vanishK} and \eqref{normderivS} (with $B$ replaced by $D_1$ and $D_2$), we obtain
\be
I = -\frac{1}{4\lambda_{-+}}e^{-2|n|s}e^{in\theta} \quad \mbox{for } n\neq 0.
\ee

We next consider $II$. From Lemma \ref{Fouriermoment}, it follows that
\begin{align*}
II=& -\frac{1}{4\lambda_{-+}}e^{-2|n|s}e^{in\theta} \\
&- \delta \sum_{m \in \mathbb{Z} \setminus \{0\}} \frac{(\varepsilon_-|nm| + \varepsilon_+nm)}{4 |n|(\varepsilon_- - \varepsilon_+)\lambda_{-+}^2}e^{-(|n|+|m|)s} \hat{h}(-n+m) e^{im\theta} + O(\delta^2),
\end{align*}
which completes the proof.
\qed

\smallskip

Consequently, we can represent the approximation of $\mathcal{A}$ in a matrix form as shown in the following result.

\begin{cor}\label{cor_asymp}
 Let $V_n$ be the subspace of $\mathcal{H}^*(\p D_2)$ spanned by $\{e^{i n \theta}, e^{-in\theta}\}$ for each $n\neq 0$.
Then the operator ${\mathcal{A}}:\mathcal{H}^*(\p D_2) \to \mathcal{H}^*(\p D_2)$ can be represented in a block matrix form
 as follows:
\begin{align}
\mathcal{A}
&=
\begin{bmatrix}
D_{1}\\
& \ddots\\
&& D_{n}\\
&&& \ddots\\
\end{bmatrix}
+\delta
\begin{bmatrix}
H_{11} & H_{12} & H_{13} & \cdots \\[1mm]
H_{21} & H_{22} & H_{23} & \cdots \\[1mm]
H_{31} & H_{32} & H_{33} & \cdots \\
\vdots & \vdots & \vdots & \ddots
\end{bmatrix}
+O(\delta^2) \nonumber\\
&=:\mathcal{A}_0 + \delta\mathcal{A}_1 + O(\delta^2).\label{pertmatrix}
\end{align}
Here, the operators $D_n: V_n \rightarrow V_n$ and $H_{nm}:V_m \rightarrow V_n$ are represented in the following matrix forms using $\{e^{in \theta},e^{-in\theta}\}$ as basis: 
\beq D_n = 
\begin{bmatrix}
\lambda^0_{-n} & 0\\
0 & \lambda^0_n\\
\end{bmatrix}
\quad \mbox{with } \lambda^0_n = -\frac{1}{4\lambda_{-+}} e^{-2|n|s}, \label{lambda_0}
\eeq
and
\beq\label{A1}
H_{nm} = -\frac{m}{4\lambda_{-+}^2}e^{-(|n|+|m|)s}
\begin{bmatrix}
2\lambda_{-+}\hat{h}(-n+m)  & \hat{h}(n+m) \\
\hat{h}(-n-m) &  2\lambda_{-+}\hat{h}(n-m) \\
\end{bmatrix},
\eeq
for $n,m \in \NN$. 

\end{cor}

The above result tells us that, in view of the eigenvalue perturbation theory \cite{kato}, we see that the eigenvalues $\lambda_n$ have the following asymptotic expansion:
\beq\label{perteigval}
\lambda_n = \lambda^0_n + \delta\lambda^1_n + O(\delta^2),
\eeq
for some coefficient $\lambda^1_n$. The eigenvalue $\lambda_n^0$ corresponds to the case of the unperturbed half-plane $\p \mathbb{H}_0$. In the case of the perturbed half-plane $\p \mathbb{H}_\delta$, the eigenvalue $\lambda_n$ is slightly different from $\lambda^0_n$.  The difference is called the shift of the eigenvalues. This shift is due to the interaction of the target particle with the plasmonic particle. 
As already explained in the previous subsection, we can measure the values of the eigenvalues $\lambda_n$. 
Since we know $\lambda^0_n$ explicitly from \eqref{lambda_0}, we can also measure the shift of the eigenvalues $\lambda_n - \lambda^0_n \approx \delta \lambda^1_n$. 
It is natural to expect that these shifts contain information on the shape of the perturbed half-plane $\p \mathbb{H}_\delta$.
In the next section, we discuss how to reconstruct the shape of $\p \mathbb{H}_\delta$ from these recovered eigenvalues.

\section{The reconstruction problem} \label{sec-inverse}

Our aim in this section is to reconstruct the shape of the small perturbation $\delta h_0$ of the half-space $\p \mathbb{H}_\delta$ using the plasmonic resonances of the small particle $\Omega$.

We assume that we measured the plasmonic eigenvalues $\lambda_n$ of the operator $\mathcal{A}$ for $n=1,2,\ldots,N$ from the far field by varying the frequency of the incident field and then picking up local peaks. Then, we can also measure the shift $\lambda_n - \lambda_n^0 \approx \delta \lambda_n^1$ as explained in Subsection  \ref{subsec:asymp}. In Proposition \ref{evencoeff}, we provide an explicit formula for the Fourier coefficients of the perturbation $\delta h$ of the transformed half-plane $\p D_1$ in terms of the shift $\delta \lambda_n^1$. This formula can be used to reconstruct the perturbation $\delta h$ in a direct way.
 Once $\delta h$ is reconstructed, the perturbation of the planar surface $\delta h_0$ can be determined by inverting the conformal map $\Phi$ and constructing $\Phi^{-1}(D_{1,\delta})$.

The following result holds.

\begin{prop}\label{evencoeff}
Let $h$ be the real-valued function in \eqref{h2}. 
We further assume that $h$ is an even function.
We have, for $n=1,2,3,\cdots$, 
\be
{\rm{Re}}\{\hat{h}(2n)\} = -\frac{2\lambda_{-+}^2(\lambda^1_{n} - \lambda^1_{-n})}{n}e^{-2ns}, \quad
{\rm{Im}}\{\hat{h}(2n)\} = 0,
\ee
and
\be
\quad \hat{h}(0) = -{\lambda_{-+}(\lambda^1_1 + \lambda^1_{-1})}e^{-2s}.
\ee
Note that, since $h$ is an even function, we have $\hat{h}(2n-1)=0$ and $\hat{h}(-n)=\hat{h}(n)$.
\end{prop}
\pf
We first derive the asymptotic expansion of the eigenvalues of the operator $\mathcal{A}$ for small $\delta$.
In what follows, by abuse of notation, we shall identify the operators $\mathcal{A},\mathcal{A}_0$ and $\mathcal{A}_1$ with their matrix representation. 
 Let us denote the standard inner product in $\ell^2(\mathbb{Z})$ by $\langle\ ,\ \rangle$. 
Let us denote by $v_n$ the eigenvector corresponding to the eigenfunction of the operator $\mathcal{A}$.
 ALso, let $\{v^0_n\}_{n\in\mathbb{Z}}$ be the standard orthonormal basis of $\ell^2(\mathbb{Z})$.
Note that $v_n^0,n\in\mathbb{Z},$ is an  eigenvector of the matrix $\mathcal{A}_0$.

We first observe that the unperturbed matrix $\mathcal{A}_0$ has degenerate eigenvalues $\lambda_n^0$ with two associated eigenvectors $v_n^0$ and $v_{-n}^{0}$. So, by applying a standard argument of the degenerate eigenvalue perturbation theory to the expansion \eqref{pertmatrix}, we see that the asymptotics of the eigenvector associated to the eigenvalue $\lambda_n$ has the following form:
\beq\label{perteigvec}
v_n = \Ga_n v_n^0 + \Gb_n v^0_{-n} + O(\delta),
\eeq
for some coefficients $\alpha_n$ and $\beta_n$.
Note that the leading order term is a linear combination of the eigenvectors $v_n^0,v_{-n}^0$ of the unperturbed matrix $\mathcal{A}_0$.

By applying \eqref{pertmatrix}, \eqref{perteigval} and \eqref{perteigvec} to
$
\Acal v_n = \Gl_n v_n,
$
 we have, up to $O(\delta^2)$,
\be
(\mathcal{A}_0 + \delta\mathcal{A}_1)(\Ga_n v_n^0 + \Gb_n v_{-n}^{0} + \delta v^1_n)=(\lambda^0_n + \delta\lambda^1_n)(\Ga_n v^0_n + \Gb_n v^0_{-n} + \delta v^1_n),
\ee
which yields the following equation by matching the first order terms with respect to $\delta$:
\beq\label{1steqn}
\Acal_0 v^1_n + \Ga_n \Acal_1 v_n^1 + \Gb_n \Acal_1 v^0_{-n} = \Gl^0_n v^1_n + \Ga_n \Gl^1_n v^0_n + \Gb_n \Gl^1_n v^0_{-n}.
\eeq
Then, by taking the inner product with $v^0_n$ and applying \eqref{pertmatrix}, we obtain
\begin{align}
\Ga_n \Gl^1_n
& = \Ga_n \langle v^0_n, \Acal_1 v^0_n \rangle + \Gb_n \langle v^0_n, \Acal_1 v^0_{-n} \rangle \nonumber\\
& = |k|\frac{\Gl^0_n}{\lambda_{-+}} \Bigl[2\Ga_n \Gl_{-+}\hat{h}(0)  + \Gb_n \hat{h}(-2n) \Bigr] \nonumber\\
& = |k|\frac{\Gl^0_n}{\lambda_{-+}} \Bigl[2\Ga_n \Gl_{-+}\hat{h}(0)  + \Gb_n \overline{\hat{h}(2n)} \Bigr],\label{acondition}
\end{align}
Similarly, taking the inner product with $v^0_{-n}$  yields
\begin{align}
\Gb_n \Gl^1_n
& = \Ga_n \langle v^0_{-n}, \Acal_1 v^0_n \rangle + \Gb_n \langle v^0_{-n}, \Acal_1 v^0_{-n} \rangle \nonumber\\
& = |n|\frac{\Gl^0_n}{\lambda_{-+}} \Bigl[\Ga_n \hat{h}(2n) + 2\Gb_n \Gl_{-+}\hat{h}(0) \Bigr],\label{bcondition}
\end{align}
where $\lambda_n^0$ is given in \eqref{lambda_0}.
Rewriting (\ref{acondition}) and (\ref{bcondition})  in a matrix form, we arrive at the following eigenvalue problem for $\lambda^1_n$:
$$
|n|\frac{\Gl^0_n}{\lambda_{-+}} \begin{bmatrix}
2\lambda_{-+} \hat{h}(0)  & \hat{h}(2n)
\\
 \hat{h}(2n) & 2\lambda_{-+} \hat{h}(0)  
\end{bmatrix}
\begin{bmatrix}
\alpha_n
\\
\beta_n
\end{bmatrix}
=
\lambda^1_{n}
\begin{bmatrix}
\alpha_n
\\
\beta_n
\end{bmatrix}.
$$
Solving this eigenvalue problem yields
\beq\label{lambda1k}
\lambda^1_{ \pm n} = |n|{\Gl_0^{|n|}}\Big(2\hat{h}(0) \pm \frac{1}{\lambda_{-+}} \hat{h}(\pm 2n)\Big).
\eeq
Note that, since $h$ is an even function, the imaginary part of every Fourier coefficient $\hat{h}(n)$ vanishes. So, we have from \eqref{hconjugate} that
$\hat{h}(-n) = \overline{\hat{h}(n)} = \hat{h}(n)$.
Therefore, \eqref{lambda1k} becomes
\beq\label{lambda1k2}
\lambda^1_{  \pm n} = |n|{\Gl^0_{|n|}}\Big(2\hat{h}(0) \pm \frac{1}{\lambda_{-+}} \hat{h}(2n)\Big).
\eeq
Then the conclusion follows from rearranging the terms with \eqref{lambda_0}. 
\qed

\begin{remark}
For the reconstruction of more general-shaped perturbations, we should be able to reconstruct the odd Fourier coefficients of $h$. These coefficients are contained in the second order terms in the asymptotic expansion for the eigenvalues of $\mathcal{A}$. So 
if the signal-to-noise ratio in the measurements is large enough, it would be possible to reconstruct the odd coefficients. The derivation of the second order term and a more complete reconstruction scheme will be considered elsewhere.
\end{remark}

\section{Numerical examples} \label{sec-numerics}

In this section, we present numerical examples to support the theoretical results. 
We assume that the measurements of the eigenvalues $\lambda_j$ 
are accurate. We generate the eigenvalues data as follows. We construct a truncated matrix for the operator $\mathcal{A}$ using formula \eqref{exact_A_mat} which is written in terms of the CGPTs. To compute these CGPTs for the transformed shape $D_{1,\delta}$, we  discretize the boundary integral in  \eqref{defCGPT} using the trapezoidal rule. We use  $300$ equispaced points for the discretization of the boundary $\p D_{1,\delta}$. By computing the eigenvalue of the resulting matrix for $\mathcal{A}$, we can obtain accurate approximations of the eigenvalues $\lambda_n$. Then, as described in Section \ref{sec-inverse}, we reconstruct the shape of $\p D_{1,\delta}$ and hence that of the perturbed half-plane $\p \mathbb{H}_\delta$.

We shall consider the case when the perturbation $\delta h$ of the planar surface is an even function. Then all the odd Fourier coefficients $\hat{h}(2n-1)$ are zero. 
From Proposition \ref{evencoeff}, we recover even Fourier coefficients of $h$. Hence, in view of \eqref{hconjugate}, we reconstruct $h$ from the following truncated Fourier series:
\be
h_{eig}(\theta) = \sum_{k=-\frac{N}{2}}^\frac{N}{2} \hat{h}(2k)e^{i2k\theta} = \hat{h}(0) + 2\sum_{k=1}^\frac{N}{2} \Re[\hat{h}(2k)]\cos(2k\theta) - 2\sum_{k=1}^\frac{N}{2} \Im[\hat{h}(2k)]\sin(2k\theta).
\ee

We compare the reconstructed images to those obtained from the complex GPTs. 
The CGPTs can be measured from the far field measurements, see \cite{book2}. It is worth mentioning that the reconstruction of the higher order CGPTs is very sensitive to measurement noise. 
We reconstruct the Fourier coefficients of $h$ from the CGPTs by using the Proposition \ref{Fouriermoment}.  

In the reconstruction of $h$ using the complex GPTs, we truncate the Fourier series of order $N$ and use \eqref{hconjugate} to arrive at
\be
h_{GPT}(\theta) = \sum_{k=-N}^N \hat{h}(k)e^{ik\theta} = \hat{h}(0) + 2\sum_{k=1}^N \Re[\hat{h}(k)]\cos(k\theta) - 2\sum_{k=1}^N \Im[\hat{h}(k)]\sin(k\theta).
\ee

In what follows, we present two numerical examples for validatation of our approach.
For their detailed implementations using MATLAB, we refer to \cite{matlab,matlab2}.

\subsection{Example 1}

We set $\delta=10^{-2}$, $d=1.04$, and $R = 0.4$. Moreover, we assume $\epsilon_- = 3$ and $\epsilon_+ = 1$. Then we get  $\lambda_{-+} = 1$ from \eqref{lambda}. Let $I$ be the union of the intervals given by $$I=[-0.7R,-0.5R] \cup [-0.3R,-0.1R]\cup[0.1R,0.3R]\cup[0.5R,0.7R].$$
We define $h_0:\RR\to\RR$ as the following even function:
 \begin{align}
 h_0(x)
 &=\ds\begin{cases}
 -1\quad &\mbox{for } x\in I,\\[5mm]
 0 \quad &\mbox{for } x \in \RR\setminus I.
 \end{cases}
 \end{align}
 Recalling $$\p \mathbb{H}_{\delta} = \{\Bx + \delta h_0(\Bx)\nu(\Bx) : \Bx \in \p \mathbb{H}_{0} \},$$
we see that the shape of the perturbed surface $\p \mathbb{H}_{\delta}$ is a corrugated plane. 
 
\begin{figure}[h]
	\centering
	\centering
    \begin{subfigure}[t]{0.45\textwidth}
	\centering
    \includegraphics[scale=0.55]{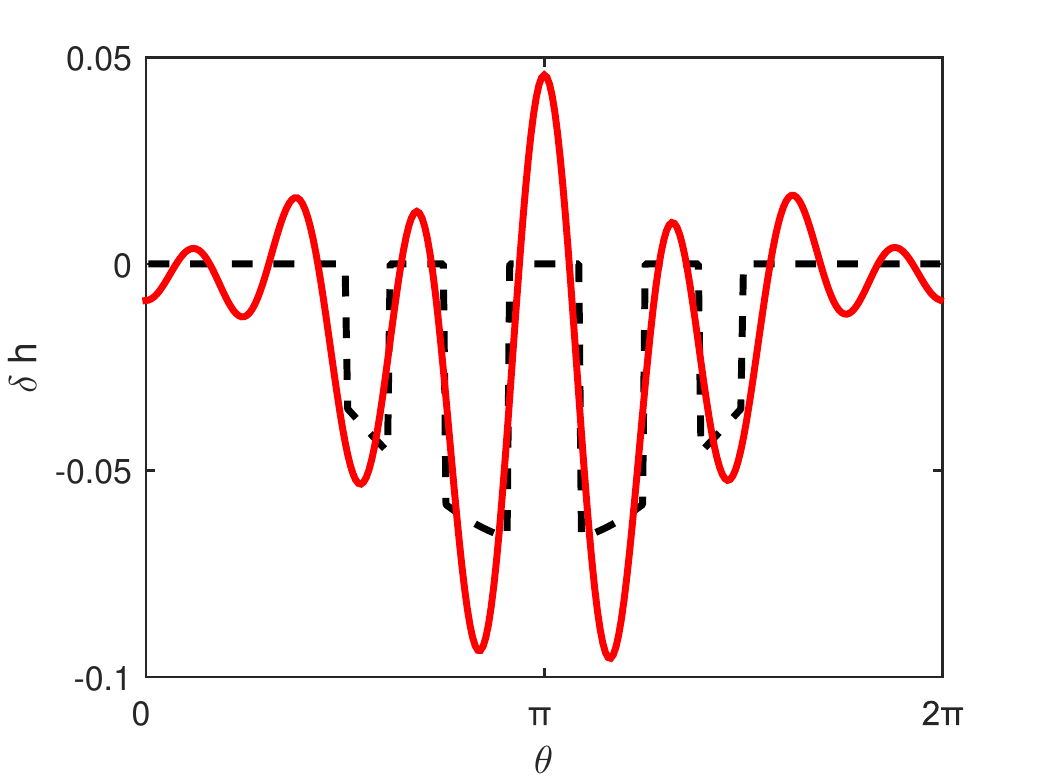}
    \caption{Reconstruction of $\delta h$ by the GPTs.}
    \label{h_GPT_ex1}
    \end{subfigure}
\quad
    \begin{subfigure}[t]{0.45\textwidth}
	\centering
    \includegraphics[scale=0.55]{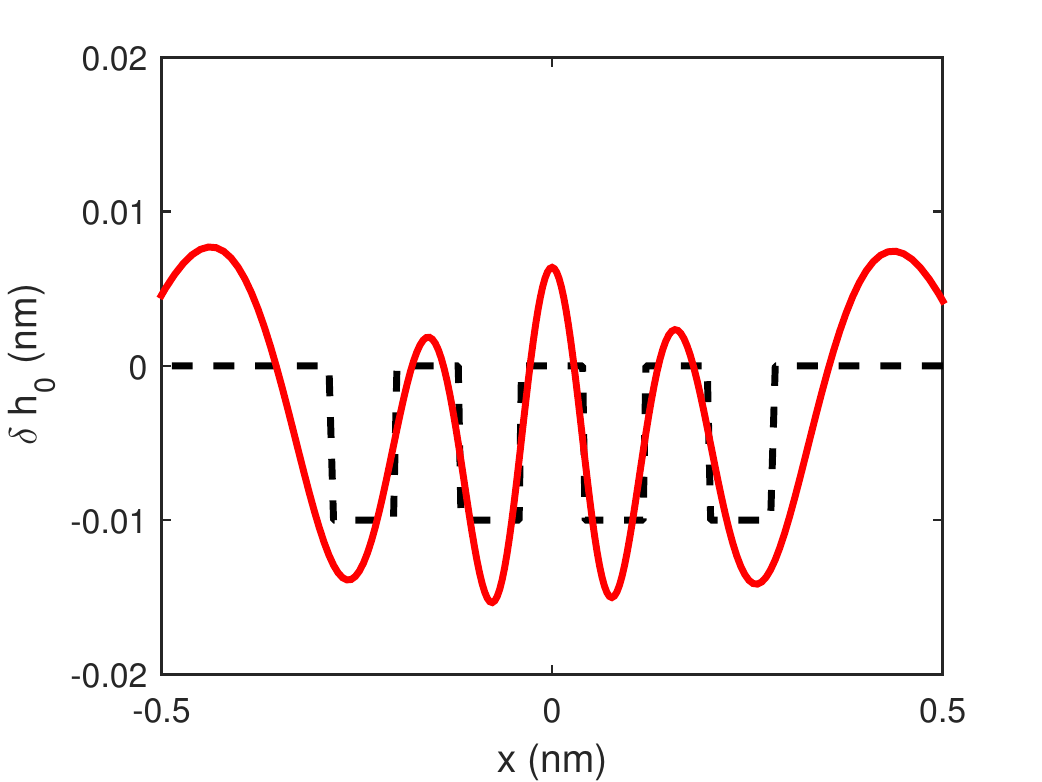}
    \caption{Reconstruction of $\delta h_0$ by the GPTs.}
    \label{h0_GPT_ex1}
    \end{subfigure}
\\
    \begin{subfigure}[t]{0.45\textwidth}
	\centering
    \includegraphics[scale=0.55]{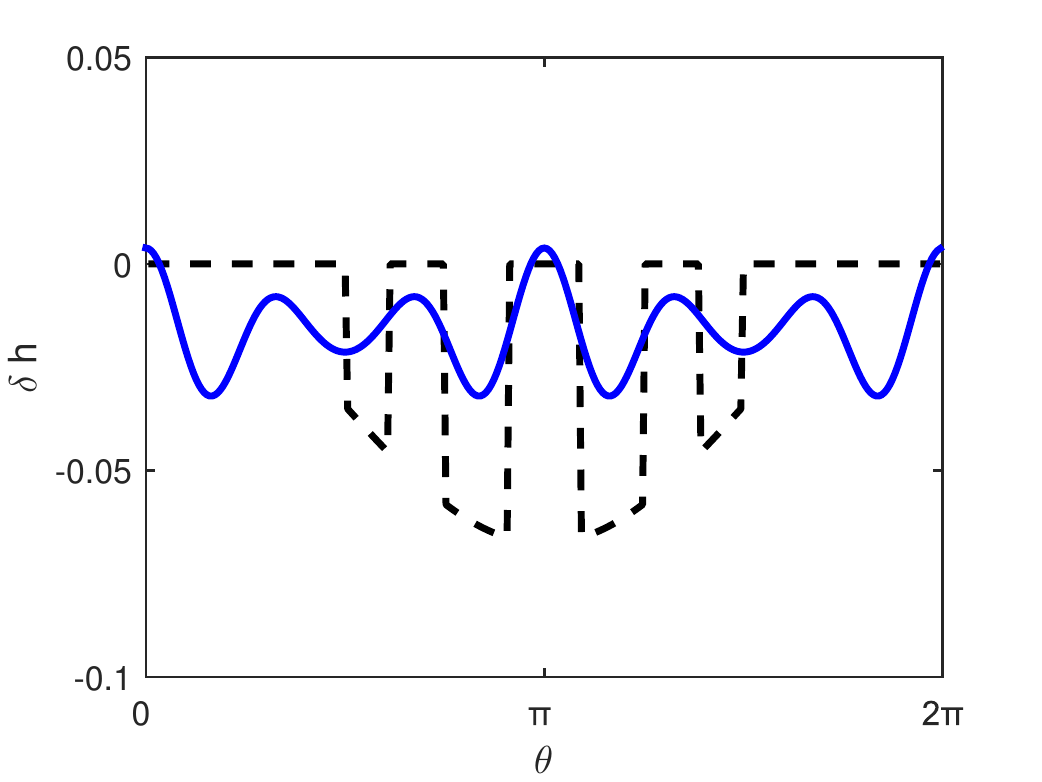}
    \caption{Reconstruction of $\delta h$ by the eigenvalues.}
    \label{h_eig_ex1}
    \end{subfigure}
\quad
    \begin{subfigure}[t]{0.45\textwidth}
	\centering
    \includegraphics[scale=0.55]{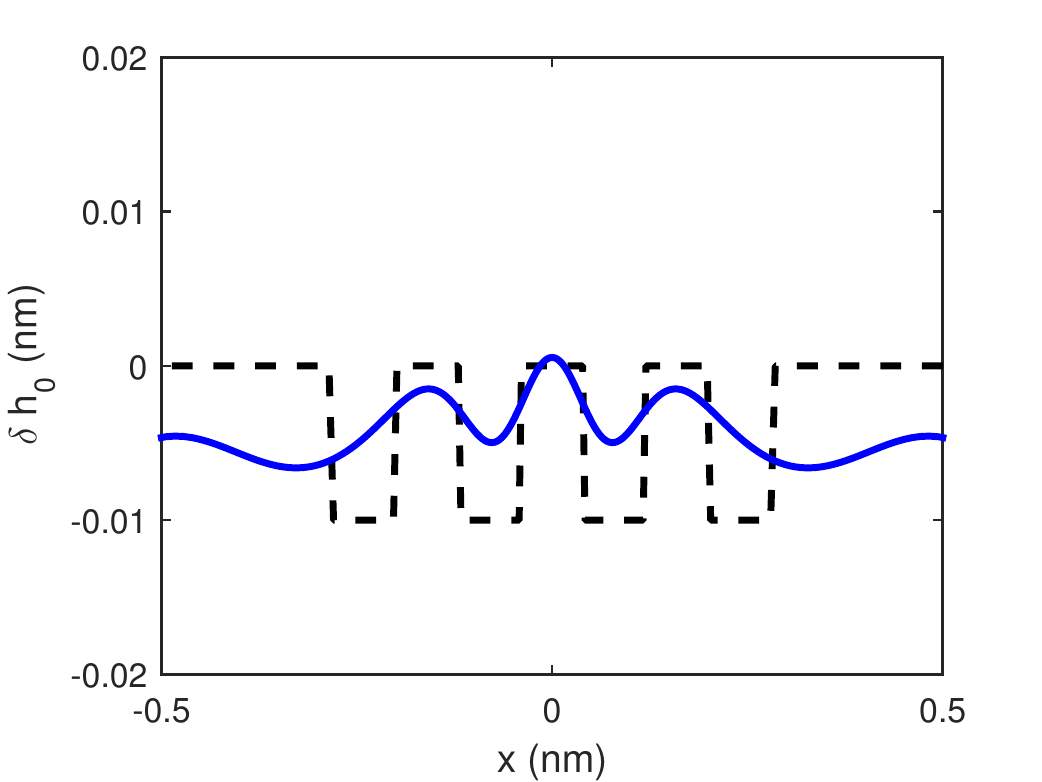}
    \caption{Reconstruction of $\delta h_0$ by the eigenvalues.}
    \label{h0_eig_ex1}
    \end{subfigure}
\caption{The dotted lines are the theoretical $\delta h:[0,2\pi]\to \mathbb{R}$ and $\delta h_0:\RR\to\RR$, the red lines are the reconstruction of $\delta h$ and $\delta h_0$ using the complex GPTs, and the blue lines are the reconstruction of $\delta h$ and $\delta h_0$ using the eigenvalue perturbations. In this case, $N=8$.}
\label{example1}
\end{figure}

Figures \ref{h_GPT_ex1} and \ref{h_eig_ex1} show the results of the reconstruction of the perturbation $\delta h_0$ based on GPTs and the shift of plasmon resonances, respectively. By transforming $\p D_{1,\delta}$ to $\p \mathbb{H}_{\delta}$ using $\Phi^{-1}$, we obtain the reconstructed $\delta h_0$ as illustrated in Figures \ref{h0_GPT_ex1} and \ref{h0_eig_ex1}.

\subsection{Example 2}

\begin{figure}[h]
	\centering
    \begin{subfigure}[t]{0.45\textwidth}
	\centering
    \includegraphics[scale=0.55]{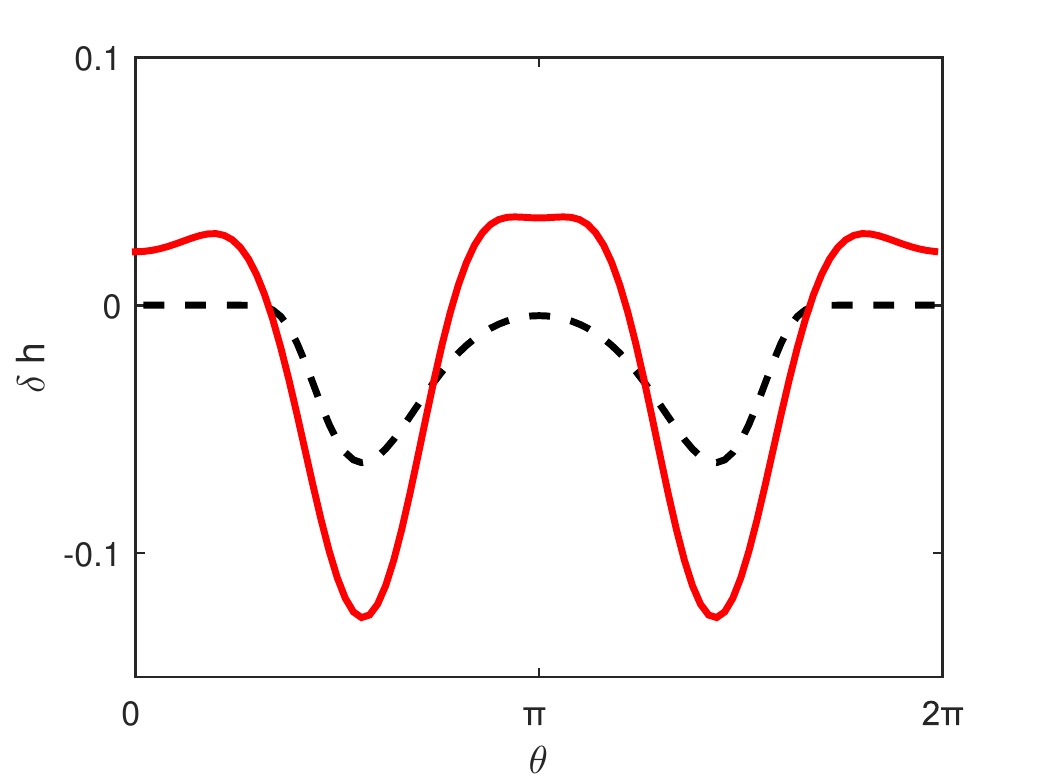}
    \caption{Reconstruction of $\delta h$ by the GPTs.}
    \label{h_GPT_ex2}
    \end{subfigure}
\quad
    \begin{subfigure}[t]{0.45\textwidth}
	\centering
    \includegraphics[scale=0.55]{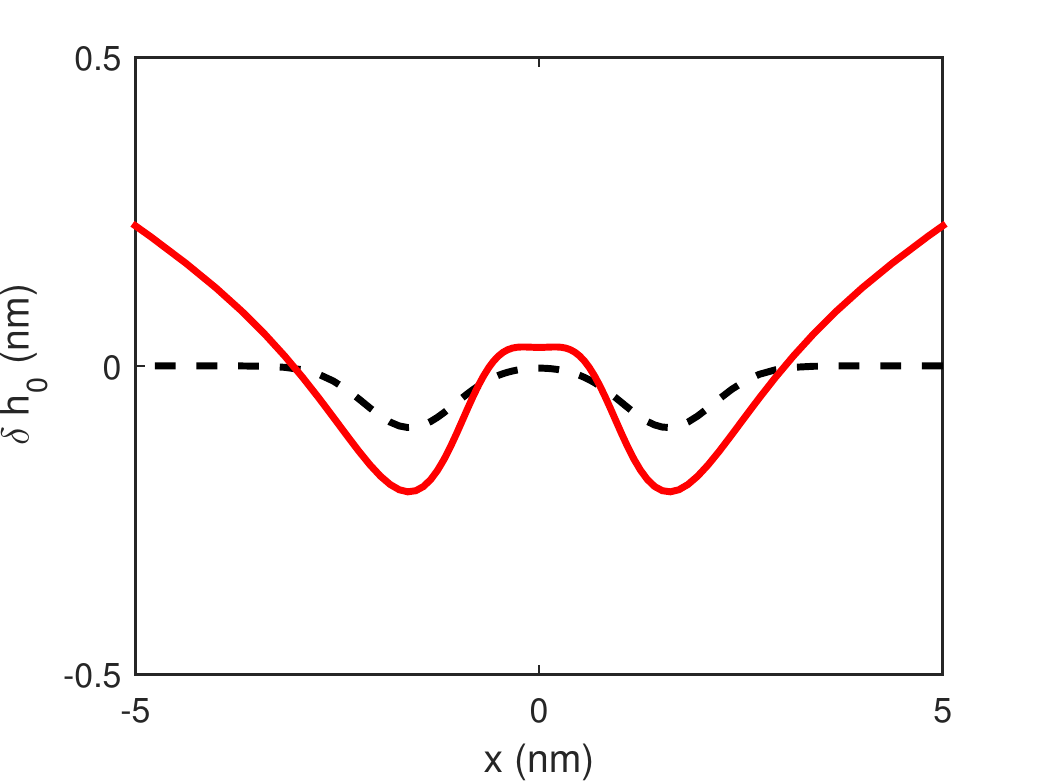}
    \caption{Reconstruction of $\delta h_0$ by the GPTs.}
    \label{h0_GPT_ex2}
    \end{subfigure}
\\
    \begin{subfigure}[t]{0.45\textwidth}
	\centering
    \includegraphics[scale=0.55]{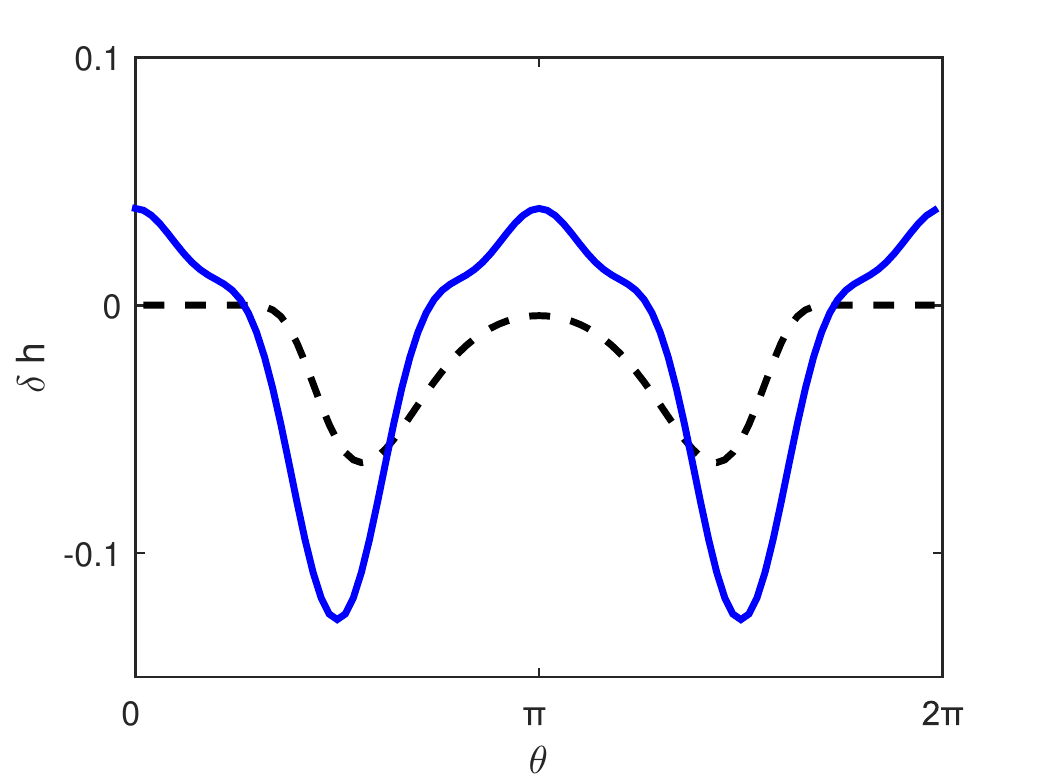}
    \caption{Reconstruction of $\delta h$ by the eigenvalues.}
    \label{h_eig_ex2}
    \end{subfigure}
\quad
    \begin{subfigure}[t]{0.45\textwidth}
	\centering
    \includegraphics[scale=0.55]{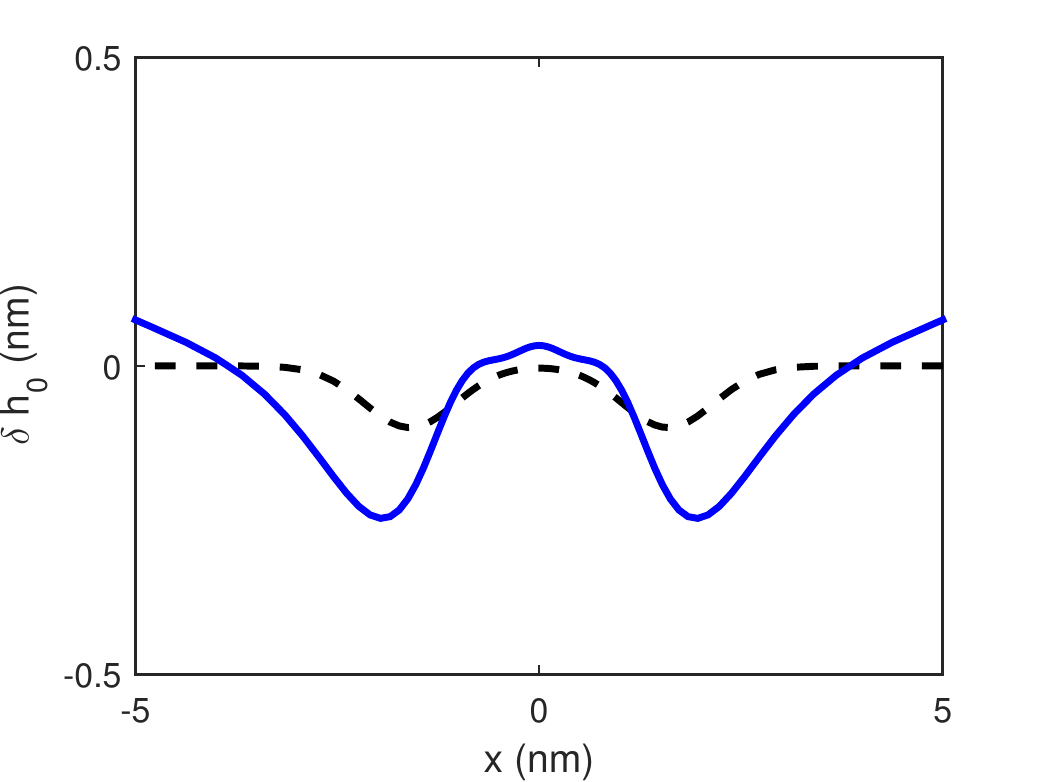}
    \caption{Reconstruction of $\delta h_0$ by the eigenvalues.}
    \label{h0_eig_ex2}
    \end{subfigure}
\caption{The dotted lines are the theoretical $\delta h:[0,2\pi]\to \mathbb{R}$ and $\delta h_0:\RR\to\RR$, the red lines are the reconstruction of $\delta h$ and $\delta h_0$ using the complex GPTs, and the blue lines are the reconstruction of $\delta h$ and $\delta h_0$ using the eigenvalue perturbations. Here, $N=6$. }
\label{example2}
\end{figure}

For the second example, we set $\delta=10^{-1}$, $d=2$ and $R = 0.8$. Moreover, we put $\epsilon_- = 3$ and $\epsilon_+ = 1$, so we have  $\lambda_{-+} = 1$ from \eqref{lambda}. We now define $h_0:\RR\to\RR$ as a smooth even function:
 \be
 h_0(x) = -\left[e^{-(\frac{x}{R}+2)^2} + e^{-(\frac{x}{R}-2)^2}\right].
 \ee
The reconstructed shapes of $\delta h$ from the CGPTs and the shifts of the eigenvalues are shown in Figures \ref{h_GPT_ex2} and \ref{h_eig_ex2}, respectively. By transforming $\p D_{1,\delta}$ to $\p \mathbb{H}_{0}$ applying $\Phi^{-1}$, we finally reconstruct $\delta h_0$ as shown in Figures \ref{h0_GPT_ex2} and \ref{h0_eig_ex2}.

In contrast with the red line in Figure \ref{h_GPT_ex1}, the blue line in Figure \ref{h_eig_ex1} shows very large fluctuations at both left and right extremities. Moreover, in comparison to the blue line in Figure \ref{h_eig_ex2}, the red line in Figure \ref{h_GPT_ex2} is much more flat at both left and right extremities. These phenomena are due to the lack of odd Fourier coefficients in $\delta h_{eig}$. Nevertheless, even Fourier coefficients of $\delta h$ contain comprehensive information of $\p \mathbb{H}_{\delta}$.

\section{Concluding remarks}

In this paper, we have introduced an original approach to recover small perturbations of a planar surface from   shifts in the resonances of the plasmonic particle-surface system. Our main idea is to design a conformal mapping which transforms the particle-surface system into a coated structure, in which the inner core corresponds to the local perturbations of the planar surface. Then we have related the perturbations of the plasmonic resonances to the Fourier coefficients of the transformed perturbations. Using these (approximate) relations, we have designed a direct (non-iterative) scheme for retrieving the perturbations of the planar surface. We have   shown that only even coefficients of the Fourier coefficients of the transformed perturbations can be reconstructed from the leading-order terms of the resonances. For large enough signal-to-noise ratio in the measurements, we may be able to recover both the odd coefficients and the complex GPTs from  second-order terms in the resonances and therefore, achieve better resolution. This would be the subject of a forthcoming work.

\end{document}